\documentclass[xypic,amscd,syntonly,amssymb]{amsart}

\usepackage{graphicx}
\usepackage[all]{xy}
\usepackage{amssymb}
\usepackage{amsmath}
\usepackage{epsfig}
\theoremstyle{plain}
\newtheorem{Thm}{Theorem}
\newtheorem{Prop}[Thm]{Proposition}
\newtheorem{Cor}[Thm]{Corollary}
\newtheorem{Lem}[Thm]{Lemma}

 \theoremstyle{definition}
\theoremstyle{remark}

\errorcontextlines=0 \numberwithin{equation}{section}

\begin{document}
 \title{Characteristic number associated to  mass linear pairs}

 \author{ ANDR\'{E}S   VI\~{N}A}
\address{Departamento de F\'{i}sica. Universidad de Oviedo.   Avda Calvo
 Sotelo.     33007 Oviedo. Spain. }
 \email{vina@uniovi.es}
\thanks{This work has been partially supported by Ministerio de Educaci\'on y
 Ciencia, grant   MAT2007-65097-C02-02}

 \maketitle
\begin{abstract}

 Let $\Delta$ be a Delzant polytope in ${\mathbb R}^n$  and ${\mathbf b}\in{\mathbb Z}^n$. Let $E$ denote the symplectic fibration over $S^2$ determined by
the pair $(\Delta,\,{\mathbf b})$. Under certain hypotheses, we
prove the equivalence between the fact that $(\Delta,\,{\mathbf
b})$ is a mass linear pair (D. McDuff, S. Tolman, {\em Polytopes
with mass linear functions. I.} Int. Math. Res. Not. IMRN 8 (2010)
1506-1574.) and the vanishing of a characteristic number of  $E$.
Denoting by ${\rm Ham}(M_{\Delta})$ the Hamiltonian group of the
symplectic manifold defined by $\Delta$, we determine loops in
 ${\rm Ham}(M_{\Delta})$ that define infinite cyclic subgroups in
 $\pi_1({\rm Ham}(M_{\Delta}))$,  when $\Delta$ satisfies any of the following
  conditions:
(i) it is the trapezium associated with a Hirzebruch surface, (ii)
it is a $\Delta_p$ bundle over $\Delta_1$, (iii) $\Delta$ is the
truncated simplex associated with the one point blow up of
${\mathbb C}P^n$.

\end{abstract}
   \smallskip

%\subjclass[2000]{ 53D50, 22E45}


%%%%%%%%%%%%%%%%%%%%%%%%%%%%%%%%%%%%%%%%%%%%%%%%%%%%%%%%%%%%%%%%%%%%%%%%%%%%%%%%%%%%%%%%%%%%%%%%%%%%%%%%%%%%%%
%%%%%%%%%%%%%%%%%%%%%%%%%%%%%%%%%%%%%%%%%%%%%%%%%%%%%%%%%%%%%%%%%%%%%%%%%%%%%%%%%%%%%%%%%%%%%%%%%%%%%%%%%%%%%%%%%%%%

\section{Introduction}

Let $(N,\Omega)$ be a closed connected symplectic $2n$-manifold.
By $\text{Ham}(N,\Omega)$, we denote the Hamiltonian group of
$(N,\Omega)$ \cite{Mc-S, lP01}. Associated with a  loop $\psi$ in
${\rm Ham}(N,\Omega)$, there exist characteristic numbers which
are invariant under deformation of $\psi$. These invariants are
defined in terms of characteristic classes of fibre bundles and
their explicit values are not easy to calculate, in general. Here,
we will consider a particular invariant $I$, whose definition we
will be recall below.
% In this article,
 By proving the
non-vanishing of $I$ for certain loops, we will deduce the
existence of infinity cyclic subgroups of $\pi_1({\rm
Ham}(N,\Omega))$, when $N$ is a toric manifold. The vanishing of
the  invariant  $I$ on particular loops in
%that are $1$-parameter subgroups
$\text{Ham}(N,\Omega)$ is related with the concept of mass linear
pair, which has been developed  in \cite{M-T3}. In this
introduction, we will state the main results of the paper and will
give a schematic exposition of the concepts involved in these
statements.

 A loop $\psi$ in
$\text{Ham}(N,\Omega)$  determines a Hamiltonian fibre bundle
$E\to S^2$ with standard fibre $N$, via the clutching
construction. Various characteristic numbers for the fibre bundle
$E$   have been defined in \cite{L-M-P}. These numbers give rise
to topological invariants of the loop $\psi$. In this article, we
will consider only the following characteristic number
 \begin{equation}\label{Ipsi}
 I({\psi}):=\int_Ec_1(VTE)\,c^n,
 \end{equation}
 where $VTE$ is the vertical tangent bundle of  $E$ and $c\in H^2(E,\,{\Bbb R})$ is the
 coupling class of the fibration $E\to S^2$ \cite{G-L-S, Mc-S}.
 $I({\psi})$ depends only on the homotopy class of the loop $\psi$.
 Moreover, the map
 \begin{equation}\label{Ihomo}
  I:\psi\in\pi_1(\text{Ham}(N,\Omega))\mapsto I({\psi})\in{\Bbb R}
   \end{equation}
 is an ${\Bbb R}$-valued group homomorphism \cite{L-M-P}.

 Our purpose is to study this characteristic number when $N$ is a
 toric manifold and $\psi$ is a $1$-parameter subgroup of ${\rm
 Ham}(N)$ defined by the toric action.
 % A toric symplectic manifold
 %is determined by its moment polytope, and
 The referred
 $1$-parameter subgroup is determined by an element ${\bf b}$ in the
 integer lattice of the Lie algebra of the corresponding torus.
  On the other hand,
 a toric symplectic manifold
is determined by its moment polytope. For a general polytope,
  a mass linear function on it is a linear
 function ``whose value on the center of mass of the polytope
 depends linearly on the positions of the  supporting hyperplanes"  \cite{M-T3}. In
 this article, we will
 %prove that a necessary and sufficient condition for
 relate the vanishing of the number $I(\psi)$ with the fact that ${\bf b}$ defines a
 mass linear function on the polytope associated with the toric
 manifold. In the following paragraphs, we provide a more detailed
 exposition of this relation.

 Let $T$ be the torus $(U(1))^n$ and $\Delta=\Delta({\bf n},k)$ the
  %Delzant
  polytope
in ${\mathfrak t}^*$ with $m$ facets defined by
\begin{equation}\label{DeltP}
\Delta({\bf n},\,k)=\bigcap_{j=1}^m\,\{ x\in {\mathfrak
t}^*\,:\,\langle x,{\bf n}_j\rangle\leq k_j \},
\end{equation}
where $k_j\in{\mathbb R}$ and the ${\bf n}_j\in {\mathfrak t}$ are
the outward conormals to the facets. The facet defined by the
equation $\langle x,\,{\bf n}_j\rangle=k_j$ will be denoted $F_j$,
and we put $\text{Cm}(\Delta)$ for the
 mass center of the polytope $\Delta$.

In \cite{M-T3}  the chamber  ${\mathcal C}_{\Delta}$ of $\Delta$
 is defined
as the set of $k'\in{\Bbb R}^m$ such that the polytope
$\Delta':=\Delta({\bf n},\,k')$ is analogous to $\Delta$; that is,
%One says that the polytopes with $m$ facets $\Delta$ and $\Delta'$
 %are analogous if
 the intersection $\cap_{j\in J} F_j$ is nonempty
iff $\cap_{j\in J} F'_j\ne \emptyset$ for any
$J\subset\{1,\dots,m\}$.
 %The polytopes in the chamber of a given
 %polytope are determined by the constants $k_j's$. Thus
When we consider only polytopes which belong to the chamber of a
fixed polytope we delete the ${\bf n}$ in the notation introduced
in (\ref{DeltP}).

 Further, McDuff and
Tolman introduced the concept of mass linear pair: Given the
polytope
  %$\tilde \Delta$
  $\Delta$ and ${\bf b}\in{\mathfrak t}$, the pair
 %$(\tilde \Delta,{\bf b})$
  $( \Delta,\, {\bf b})$ is mass linear if the map
\begin{equation}\label{kmapsto}
k\in{\mathcal C}_{\Delta}\mapsto \langle
\text{Cm}(\Delta(k)),\,{\bf b}\rangle\in{\Bbb R}
 \end{equation}
 is linear. That is,
 \begin{equation}\label{Rjkj}
 \langle {\rm Cm}(\Delta(k),\,{\bf b}\rangle=\sum_j R_jk_j+C,
 \end{equation}
where $R_j$
 and $C$ are constant.

Let us assume that $\Delta$ is a Delzant polytope \cite{Del}. We
shall denote by $(M_{\Delta},\,\omega_{\Delta},\,\mu_{\Delta})$
the toric manifold determined by $\Delta$
 %\cite{Gui}, and with
 ($\mu_{\Delta}:M\to{\mathfrak t}^*$  being the corresponding moment map).
 %(see \cite{Gui}).
 Given ${\bf b}$, an element in the integer lattice of ${\mathfrak
t}$, we shall write $\psi_{\bf b}$ for the loop of Hamiltonian
diffeomorphisms of $(M_{\Delta},\,\omega_{\Delta})$ defined by
${\bf b}$ through the toric action.
 We
will let $I(\Delta;\,{\bf b})$ for the characteristic number
$I(\psi_{\bf b})$.
 When we consider only
polytopes
 in the chamber of a given polytope, we will write $I(k;\,{\bf b})$
 instead of $I(\Delta(k);\,{\bf b})$ for $k$ in this chamber.

The group $G$ of the translations defined by the elements of
${\mathfrak t}^*$ acts freely on ${\mathcal C}_{\Delta}$. We put
$r:=m-n$ for the dimension of the quotient ${\mathcal
C}_{\Delta}/G.$ Thus, $r$ is the number of effective parameters
which characterize the polytopes in ${\mathcal C}_{\Delta}$
considered as ``physical bodies".

We will prove the following theorem:

\begin{Thm}\label{FinalThm}
Let   $(\Delta,\,{\bf b})$ be a pair consisting of a Delzant
polytope in ${\mathfrak t}^*$ and an element in the integer
lattice of ${\mathfrak t}$. If $r\leq 2$, the following statements
are equivalent

(a)  $I(k;\,{\bf b})=0,$ for all $k\in{\mathcal C}_{\Delta}$.

(b)
% $(\Delta,\,{\bf b})$ is a  mass linear pair.
  $(\Delta,\,{\bf b})$ is a  mass linear pair as in
(\ref{Rjkj}), with $\sum_jR_j=0$.
 \end{Thm}

In \cite{Vipp}, by direct computation, we proved the
  equivalence  between the vanishing of $I(k;\,{\bf b})$ on
  ${\mathcal C}_{\Delta}$ and the fact that $(\Delta,\,{\bf b})$ is a mass
  linear pair,  when $\Delta$ satisfies any of the following
  conditions:

  (i) it is the trapezium associated with a Hirzebruch surface,

(ii) it is a $\Delta_p$ bundle over $\Delta_1$ \cite{M-T3},

(iii)   $\Delta$ is the truncated simplex associated with the one
point blow up of ${\mathbb C}P^n$.

\smallskip

 On the other hand, when $\Delta$ is
any of these polytopes (i)-(iii), the number $r$ is equal to $2$;
thus, from Theorem \ref{FinalThm} and the result of  \cite{Vipp},
it follows that condition $\sum_j R_j=0$ is satisfied by
 all the mass linear pairs $(\Delta,\,{\bf b})$. This fact can
 also be proved  by
  direct calculation (Propositions \ref{PropHirz}, \ref{SumRjBundle} and  \ref{Propsumrblow}).
  So,
Theorem \ref{FinalThm}, together with these Propositions,
generalize the
 result proved in    \cite{Vipp}.

Although the homotopy type of the Hamiltonian groups ${\rm
Ham}(N,\Omega)$ is known only for some symplectic manifolds
\cite{McD},
 the invariant $I$ allows us to
identify nontrivial elements in $\pi_1({\rm Ham}(N,\Omega))$. As
$I$ is a {\em group homomorphism}, from Theorem \ref{FinalThm}, we
deduce that a sufficient condition for $\psi_{\bf b}$ to generate
an infinite cyclic subgroup in $\pi_1( {\rm
Ham}(M_{\Delta},\,\omega_{\Delta}))$ is that the above condition
(b) does not hold for $(\Delta,\,{\bf b})$.
%, assumed that $r\leq 2$.
More precisely, we have the following consequence of Theorem
\ref{FinalThm}:

\begin{Thm}\label{CorThm1}
 Given the Delzant polytope $\Delta$ and ${\bf b}$ an element in the
 integer lattice of ${\mathfrak t}$. If $r\leq 2$ and $(\Delta,\,{\bf
 b})$ is not mass linear, then $\psi_{\bf b}$ generates an infinite cyclic subgroup in $\pi_1(
{\rm Ham}(M_{\Delta(k)},\,\omega_{\Delta(k)}))$, for all
$k\in{\mathcal C}_{\Delta}$.
 \end{Thm}

In the proof of Theorem \ref{FinalThm},  a formula for
 the characteristic number $I({\psi_{\bf b}})$ obtained in
\cite{V2} plays a crucial role. This formula gives $I({\psi_{\bf
b}})$
 %the value of the characteristic number
  in terms of the integrals, on the facets of
the polytope, of the normalized Hamiltonian function corresponding
to the loop $\psi_{\bf b}$ (see (\ref{Ibr})). From this expression
for $I({\psi_{\bf b}})$, we will deduce a relation between the
directional derivative of map (\ref{kmapsto}) along the vector
$(1,\dots,1)$ of ${\mathbb R}^m$, the Euclidean volume of
$\Delta(k)$ and $I(k;\,{\bf b})$ (see (\ref{util})).
 From this relation, it is easy to complete the proof of Theorem
 \ref{FinalThm}.

This article is organized as follows: In Section 2,  we   study
the characteristic number $I(k;\,{\bf b})$, when $(\Delta,\,{\bf
b})$ is a linear
 pair and $k$ varies in the chamber of $\Delta$;
 we    prove that $I(k;\,{\bf b})$ is a homogeneous polynomial of the $k_j$
 (Proposition \ref{polyno}).

In Section 3, we prove Theorem \ref{FinalThm}.
 In Proposition \ref{PropSufCon}, a sufficient geometric condition for the
 Delzant polytope $\Delta$ to admit a mass linear pair
 $(\Delta,\,{\bf b})$   is given.
  For a Delzant polytope $\Delta$, Proposition \ref{PropNecCond}
  gives a necessary condition for the vanishing of $I(k;\,{\bf b})$   on ${\mathcal C}_{\Delta}$.
 We also
 %give an
 %equivalent condition to the vanishing of
 express $\sum_j R_j$ in terms of
the displacement of the center of mass ${\rm Cm}(\Delta(k))$
produced by the change $k_j\to k_j+1$ (Proposition
\ref{PropsumR}).

 Section 4  concerns   the   form which Theorem \ref{CorThm1} adopts, when $\Delta$
 is a Delzant polytope of the particular types (i)-(iii)
 % considered in \cite{Vipp} and that have been
mentioned above (see Corollary \ref{CoroHirz}, Theorems
\ref{pi1fifration} and  \ref{Propblowup}).
 We also prove that, in these particular cases, if $(\Delta,{\bf
b})$ is a mass linear pair, then $\sum_jR_j=0$.
 %(Propositions \ref{PropHirz},  \ref{SumRjBundle} and \ref{Propsumrblow}).

%%%%%%%%%%%%%%%%%%%%%%%%%%%%%%%%%%%%%%%%%%%%%%%%%%%%%%%%%%%%%%%%%%%%%%%%%%%%%%%%%%%%%%%%%%%%%%%%%%%%%%%%%%%%%%%%%%%%%%%%%%%%%
%%%%%%%%%%%%%%%%%%%%%%%%%%%%%%%%%%%%%%%%%%%%%%%%%%%%%%%%%%%%%%%%%%%%%%%%%%%%%%%%%%%%%%%%%%%%%%%%%%%%%%%%%%%%%%%%%%%%%%

\section{A characteristic number}\label{SecChar}

Let us suppose that the polytope $\Delta$  defined in
(\ref{DeltP}) is a Delzant polytope
 in ${\mathfrak t}^*$.  Following \cite{Gui}, we recall some points of the construction of
$(M_{\Delta},\,\omega_{\Delta},\,\mu_{\Delta})$ from the polytope
$\Delta$. We put   $\tilde T:=(S^1)^{m-n}$. The ${\bf n}_i$
determine weights $w_j\in\tilde{\mathfrak t}^*$, $j=1,\dots,m$ for
a $\tilde T$-action on ${\Bbb C}^m$. Then
 %the map $J$ is a
moment map for this action is
$$J:z\in{\Bbb C}^m\mapsto J(z)=\pi\sum_{j=1}^m|z_j|^2w_j\in \tilde{\mathfrak
t}^*.$$
 The $k_i$ define a regular value $\sigma$ for $J$, and the
 manifold $M_{\Delta}$ is the following orbit space
 \begin{equation}\label{DefiM}
 M_{\Delta}=\big\{z\in {\Bbb C}^m\,:\,\pi\sum_{j=1}
^m|z_{j}|^2w_{j}=\sigma \big\}/{\tilde T},
 \end{equation}
 where the relation defined by ${\tilde T}$ is
 \begin{equation}\label{DefiAc}(z_j)\simeq(z'_j)\;\;\text{iff
there is}\;\,{\bf y}\in \tilde {\mathfrak t}\;\, \text{such
that}\;\, z'_j=z_j\,e^{2\pi i\langle w_j,{\bf y} \rangle}\;
\,\text{for}\;\, j=1,\dots, m.
 \end{equation}

Identifying $\tilde{\mathfrak t}^*$ with ${\Bbb R}^r$,
$\sigma=(\sigma_1,\dots,\sigma_r)$ and each $\sigma_a$ is a linear
combination of the $k_j$.

 Given a facet $F$ of $\Delta$, we choose a vertex $p$ of $F$. After a possible change in numeration
 of the facets, we can assume that $F_1,\dots,F_n$ intersect at $p$. In this numeration $F=F_j$, for some
 $j\in\{1,\dots,n\}$.

  If we write $z_a=\rho_a e^{i\theta_a}$, then
the symplectic form can be written on $\{[z]\in M\,:\, z_a\ne 0,\,
\forall a \}$
 \begin{equation}\label{omeg}
  \omega_{\Delta}=(1/2)\sum_{i=1}^n d\rho_i^2\wedge
d\varphi_i,
 \end{equation}
  with $\varphi_i$ an angular variable, linear combination of the $\theta_a$.

The action of $T=(S^1)^n$ on $M_{\Delta}$
$$(\alpha_1,\dots,\alpha_n)[z_1,\dots,
z_m]:=[\alpha_1z_1,\dots,\alpha_nz_n,z_{n+1},\dots,z_m]$$
 endows   $M_{\Delta}$ with a structure of  toric manifold.
Identifying ${\mathfrak t}^*$ with ${\mathbb R}^{n}$,
 the moment map
 $\mu_{\Delta}:M_{\Delta}\to{\mathfrak t}={\Bbb R}^n$ is defined by
 \begin{equation}\label{dd}
  \mu_{\Delta}([z])=\pi(\rho^2_1,\dots,\rho_n^2)+(d_1,\dots,d_n),
 \end{equation}
where the constants $d_i$  are linear combinations of the $k_j$
and
 \begin{equation}\label{immu}
  \text{im}\,\mu_{\Delta}=\Delta.
 \end{equation}
The facet $F=F_j$ of $\Delta$ is the image by $\mu_{\Delta}$ of
the submanifold
 $$D_j=\{ [z_1,\dots,z_m]\in M_{\Delta} \,|\,
z_j=0\}.$$

We write $x_i:=\pi\rho_i^2$, then
 \begin{equation}\label{Intpolytope1}
 \int_{M_{\Delta}}(\omega_{\Delta})^n=
n! \int_{\Delta} dx_1\dots dx_n.
\end{equation}

Let ${\bf b}$ be an element in the integer lattice of ${\mathfrak
t}$. The normalized Hamiltonian of the   circle action generated
by ${\bf b}$ is the function $f$ determined by,
 $$f=\langle \mu_{\Delta},\,{\bf b}\rangle+ \text{constant}\;\;\;\text{and}\;\;\;
\int_{M_{\Delta}}f\,(\omega_{\Delta})^n=0.$$
 That is, $f=\langle \mu_{\Delta},\,{\bf b}\rangle-\langle{\rm
 Cm}(\Delta),\,{\bf b}\rangle,$ where
 \begin{equation}\label{Cm}
\langle \text{Cm}(\Delta),\,{\bf b}\rangle= \frac{\int_M
\langle\mu_{\Delta},\,{\bf b}\rangle\,(\omega_{\Delta})^n}{\int_M
(\omega_{\Delta})^n}.
\end{equation}
Moreover,
\begin{equation}\label{Intpolytope2}
  \int_{M_{\Delta}}\langle\mu_{\Delta},\,{\bf b}\rangle(\omega_{\Delta})^n=
 n!\int_{\Delta}\sum_{i=1}^n b_ix_i\, dx_1\dots dx_n.
 \end{equation}

An expression    for the value of the invariant $I({\psi_{\bf
b}})$ in terms of integrals of the Hamiltonian function has been
obtained in Section 4 of \cite{V2} (see also \cite{V} and
\cite{Sh})
\begin{equation}\label{Ibr}
I(\Delta;\,{\bf b}):=I({\psi_{\bf b}})=-n\sum_{F\,{\rm
facet}}N(F),
\end{equation}
 where the contribution $N(F)$ of the above facet $F=F_j$ (with $j=1,\dots, n$) is
    \begin{align}\label{N(F)}
    N_j:&=N(F)=(n-1)!\int_{F_j}f \,dx_1\dots \hat{dx_j}\dots dx_n \\
    \notag
   & =(n-1)!\Big( \int_{F_j}\langle\mu_{\Delta},\,{\bf b}\rangle\,dx_1\dots \hat{dx_j}\dots
    dx_n     -\langle {\rm Cm}(\Delta),\,{\bf b}\rangle\int_{F_j}dx_1\dots \hat{dx_j}\dots
    dx_n\Big),
\end{align}
with $dx_1\dots \hat{dx_j}\dots
    dx_n:=dx_1\dots dx_{j-1}dx_{j+1}\dots
    dx_n.$

Given $\Delta=\Delta({\bf n},\,k)$, we consider the polytope
 $\Delta'=\Delta({\bf n},\,k')$ obtained from
$\Delta$ by the translation defined by a vector $a$ of ${\mathfrak
t}^*$. As we said, we write $I(k;\,{\bf b})$ and $I(k';\,{\bf b})$
for the corresponding characteristic numbers. According to the
construction of the respective toric manifolds,
$$M_{\Delta'}=M_{\Delta},\;\;\;
\omega_{\Delta'}=\omega_{\Delta},\;\;\;
\mu_{\Delta'}=\mu_{\Delta}+a.$$
  But the {\em normalized}
Hamiltonians $f$ and $f'$ corresponding to the action of ${\bf b}$
on $M_{\Delta}$ and $M_{\Delta'}$ are equal. Thus, it follows from
(\ref{Ibr}) that $I(k;\,{\bf b})=I(k';\,{\bf b})$. More precisely,
we have the evident proposition:

\begin{Prop}\label{Propinv}
If $a$ is an arbitrary vector of ${\mathfrak t}^*$, then
$I(k;\,{\bf b})=I(k';\,{\bf b})$, for $k'_j=k_j+\langle a,{\bf
n}_j\rangle$, $\;j=1,\dots, m$.
\end{Prop}

 By Proposition \ref{Propinv},  we can assume that all $d_j$ in (\ref{dd}) are
 zero for the determination of $I(k;\,{\bf b})$.

The following lemma is elementary:
\begin{Lem}\label{Lemmaint}
If
 $$S_n(\tau):=\Big\{(x_1,\dots,x_n)\in{\mathbb
R}^n\,\Big|\,\sum_{i=1}^nx_i\leq\tau,\;\;\;0\leq x_j, \;\forall
j\Big\},$$ then
$$\int_{S_n(\tau)}f(x_1,\dots,x_n)\,dx_1\dots dx_n=\begin{cases}

 \frac{\tau^{n}}{n!}\, ,\; \text{if}\;\;\; f=1 \\

                                         \\

c\frac{\tau^{n+c}}{(n+c)!}\, ,\; \text{if}\;\;\; f=x_i^c,\; c=1,2 \\
              \\
\frac{\tau^{n+2}}{(n+2)!}\, ,\; \text{if}\;\;\; f=x_i x_j,\;  i\ne
j.
\end{cases}$$

\end{Lem}

More general, if $c_1,\dots,c_n\in{\mathbb R}_{>0}$, we put
$$S_n(c,\tau):=\Big\{(x_1,\dots,x_n)\in{\mathbb
R}^n\,\Big|\,\sum_{i=1}^n c_ix_i\leq\tau,\;\;\;0\leq x_j,
\;\forall j\Big\},$$
 then
 \begin{equation}\label{Intsnctau}
  \int_{S_n(c,\tau)}dx_1\dots
 dx_n=\frac{1}{n!}\prod_{i=1}^n\frac{\tau}{c_i},\;\;\; \; \int_{S_n(c,\tau)}x_j\,dx_1\dots
 dx_n=\frac{1}{(n+1)!}\frac{\tau}{c_j}\prod_{i=1}^n\frac{\tau}{c_i}
  \end{equation}
Thus, in the particular case  that $\Delta=S_n(c,\tau)$, the
integral $\int_{M_{\Delta}}(\omega_{\Delta})^n$ is a monomial of
degree $n$ in $\tau$, and $\int_{M_{\Delta}}\langle
\mu_{\Delta},\,{\bf b}\rangle(\omega_{\Delta})^n$ is a monomial of
degree $n+1$.

\smallskip

We return to the general case in which  $\Delta$ is the polytope
defined in (\ref{DeltP}). Its
 vertices   are the solutions to
 \begin{equation}\label{verticed(Delta)}
 \langle x,\,{\bf
n}_{j_a}\rangle=k_{j_a},\;\, a=1,\dots,n;
 \end{equation}
 hence, the
coordinates of any vertex of $\Delta$ are linear combinations of
the $k_j$.

A hyperplane in ${\mathbb R}^n$ through a vertex
$(x^0_1,\dots,x_n^0)$ of $\Delta$ is given by an equation of the
form
 \begin{equation}\label{langlex}
  \langle x,\,{\bf n}\rangle =\langle x^0,\,{\bf
n}\rangle=:\kappa.
 \end{equation}
 Thus, the independent term $\kappa$ is a linear combination (l. c.) of the $k_j$. Moreover, the
 coordinates of the
 common point of $n$ hyperplanes
 \begin{equation}\label{hyperpl.c.}
 \langle x,\,{\bf \tilde n}_i\rangle =\kappa_i,
  \end{equation}
  with $\kappa_i$
 l. c. of the $k_j$   are also
 l. c. of the $k_j$.

 By drawing  hyperplanes   through vertices of $\Delta$ (or more generally, through points which are the
 intersection of $n$ hyperplanes as (\ref{hyperpl.c.})),   we can obtain
 a family $\{\,_{\beta}S\}$ of subsets of $\Delta$ such that:

 a) Each $_{\beta}S$ is the transformed of a simplex
 $S_n(b,\tau)$ by an element of the
 group of Euclidean motions in ${\mathbb R}^n$.

 b) For $\alpha\ne\beta$, $\, \,_{\alpha}S\cap\,_{\beta}S$  is a
 subset of the border of $_{\alpha}S$.

 c) $\bigcup_{\beta}\,_{\beta}S=\Delta$.

Thus, by construction, each facet of $_{\beta}S$ is contained in a
hyperplane $\pi$ of the form $\langle x,\,{\bf n}\rangle =
\kappa$,
 with $\kappa$ l. c. of the $k_j$.

On the other hand, the hyperplane $\pi$  is transformed by an
element of ${\rm SO}(n)$ in an hyperplane $\langle x,\,{\bf
n'}\rangle=\kappa$.  If ${\mathcal T}$ is a translation in
${\mathbb R}^n$ which applies $S_n(b,\tau)$ onto $_{\beta}S$, then
this transformation maps $(0,\dots, 0)$ in a vertex
$a=(a_1,\dots,a_n)$ of $_{\beta}S$. So, the translation ${\mathcal
T}$ transforms $\pi$ in $\langle x,\,{\bf n}\rangle=\kappa+\langle
a,\,{\bf n}\rangle=:\kappa'$. As each $a_j$ is a l. c.  of the
$k_j$, so is $\kappa'$. Hence, any element of the group of
Euclidean motions in ${\mathbb R}^n$ which maps $S_n(b,\tau)$ onto
$_{\beta}S$ transforms the hyperplane $\pi$
 %through a vertex of $\Delta$ in an hyperplane
\begin{equation}\label{hyperplanes}
\langle x,\,{\bf n'}\rangle=\kappa',
 \end{equation}
 with $\kappa'$
a l. c. of the $k_j$.

 Let assume that  $(R{\mathcal T}_a)(S(b,\tau))= {}_{\beta}S$, with
 $R\in{\rm SO}(n)$ and ${\mathcal T}_a$  the translation  defined by $a$.
  Then the oblique facet of
$S(b,\tau)$,
 contained in the hyperplane $\sum b_ix_i=\tau$, is the image by
 $T_{-a}R^{-1}$ of a facet of $_{\beta}S$, which in turn is contained in a
 hyperplane of equation (\ref{hyperplanes}) ($\kappa'$ being
a l. c. of the $k_j$). The argument of the
 preceding paragraph applied to $R^{-1}$ and ${\mathcal T}_{-a}$
 proves that
$\tau$ is a l. c. of the $k_j$. Hence,   by (\ref{Intsnctau}) the
integral
 $$\int_{_\beta S}dx_1\dots dx_n=\int_{S_{n}(b,\tau)}dx_1\dots
 dx_n$$
 is a monomial of degree $n$ of a l. c. of the $k_j$.
 Thus,
 \begin{equation}\label{IntHomog}
 \int_{M} (\omega_{\Delta})^n=\sum_{\beta}\int_{_\beta
 S}dx_1\dots dx_n,
  \end{equation}
is a homogeneous polynomial of degree $n$ of the $k_j$.

Similarly,
 \begin{equation}\label{IntHomog2}
 \int_{M_{\Delta}}\langle\mu_{\Delta},\,{\bf
b}\rangle(\omega_{\Delta})^n
\end{equation}
 is a homogeneous polynomial of degree $n+1$ of the $k_j$.
 Analogous results hold for
$$\int_{F_j}dx_1\dots\hat{dx_j}\dots dx_n,\;\; \hbox{and}\;\;
\int_{F_j}\langle \mu_{\Delta},\,{\bf b}\rangle\,
dx_1\dots\hat{dx_j}\dots dx_n.$$

 From   formulas (\ref{Intpolytope1})-(\ref{N(F)}) together with the
 preceding argument, it follows
 % together with Lemma \ref{Leml_ilinear}
   the following
 proposition:
 \begin{Prop}\label{HomPr}
 Given a Delzant polytope $ \Delta$, if ${\bf b}$ belongs to
 the integer lattice of ${\mathfrak t}$, then
 $I(k;\,{\bf b})$
is a rational function of the $k_j$, for $k\in{\mathcal
C}_{\Delta}$.
\end{Prop}

 Analogously, we have

\begin{Prop}\label{polyno}
If $(\Delta,\,{\bf b})$ is mass linear pair, then $I(k;\,{\bf b})$
is a homogeneous polynomial in the $k_j$ of degree $n$, when
$k\in{\mathcal C}_{\Delta}$.
\end{Prop}

We will use the following simple lemma in the proof of Theorem
\ref{FinalThm}.
\begin{Lem}\label{LemCm}
If $\Hat k_j=sk_j$ for $j=1,\dots,m$, with $s\in{\mathbb R}$, then
${\rm Cm}(\Delta({\bf n},\,\Hat k))=s\,{\rm Cm}(\Delta({\bf n},\,
k))$.
\end{Lem}
{\it Proof.} The vertices of $\Delta ({\bf n},\,k)$ are the
solutions of (\ref{verticed(Delta)}) and the vertices of $\Delta
({\bf n},\,\Hat k)$ are the solutions of $\langle x,\,{\bf
n}_{j_a}\rangle=sk_{j_a}$, with $a=1,\dots,n$. Thus, the vertices
of $\Delta ({\bf n},\,\Hat k)$ are those of $\Delta ({\bf n},\,k)$
multiplied by $s$. \qed

The Lemma also follows from the fact that (\ref{IntHomog}) and
(\ref{IntHomog2}) are homogeneous polynomials of degree $n$ and
$n+1$, respectively.

\smallskip

%%%%%%%%%%%%%%%%%%%%%%%%%%%%%%%%%%%%%%%%%%%%%%%%%%%%%%%%%%%%%%%%%%%%%%%%%%%%%%%%%%%%%%%%%%%%%%%%%%%%%%%%%%%%%%%%%%
%%%%%%%%%%%%%%%%%%%%%%%%%%%%%%%%%%%%%%%%%%%%%%%%%%%%%%%%%%%%%%%%%%%%%%%%%%%%%%%%%%%%%%%%%%%%%%%%%%%%%%%%%%%%%%%%%%%%%%%

\section{Proof of Theorem \ref{FinalThm}}

Let  us assume that the polytope $\Delta$ defined by (\ref{DeltP})
is  Delzant
% polytope   in ${\mathfrak t}^*$
 %defined in (\ref{DeltP})
  and let $k$ be an element of ${\mathcal C}_{\Delta}$. We denote by
$M_{(k)}$, $\,\omega_{(k)}$ and $\mu_{(k)}$, the manifold, the
symplectic structure and the moment map (resp.) determined by
$\Delta(k)$. The facets of $\Delta(k)$ will be denoted by by
$F_{(k)j}$.
% and we set
%$D_{(k)j}:=\mu_{(k)}^{-1}(F_{(k)j})$.

Let ${\bf b}$ be an element in the integer lattice of ${\mathfrak
t}$. We put
$$A_{(k)}:=\int_{M_{(k)}}\langle\mu_{(k)},\,{\bf b}\rangle(\omega_{(k)})^n,\;\;\;\;\;
B_{(k)}:=\int_{M_{(k)}}(\omega_{(k)})^n.$$

By (\ref{Intpolytope1}), $\frac{1}{n!}B_{(k)}$ is the Euclidean
volume of the polytope $\Delta(k)$. Given a facet $F_{(k)j}$, we
can assume that $j\in\{1,\dots,n\}$ (see  third paragraph of
Section \ref{SecChar}). So, $F_{(k)j}$ is defined by the equation
$x_j=0$. If we make an infinitesimal variation of the facet
$F_{(k)j}$, by means of the translation defined by $k_j\to
k_j+\epsilon$ (keeping unchanged the other $k_i$), then the volume
of $\Delta(k)$ changes according to
$$\frac{1}{n!}B_{(k)}\longrightarrow
\frac{1}{n!}B_{(k)}+\epsilon \int_{F_{(k)j}}dx_1\dots\hat
{dx_j}\dots dx_n+ O({\epsilon}^{2}).$$ We write $dX^j$ for
$dx_1\dots\hat {dx_j}\dots dx_n$.
% since
%$$\frac{1}{(n-1)!}\int_{D_{(k)j}}(\omega_{(k)})^{n-1}$$
%is the Euclidean volume of $F_{(k)j}$.
 Thus,
% From (\ref{Intpolytope}) it follows,
 $$\frac{\partial B_{(k)}}{\partial
 k_j}=n!\int_{F_{(k)j}}dX^j,\;\;\;\;\;
  \frac{\partial A_{(k)}}{\partial
 k_j}=n!\int_{F_{(k)j}} \langle\mu_{(k)},\,{\bf b} \rangle\, dX^j. $$
 So, by (\ref{Cm}),
$$\frac{\partial}{\partial k_j}\langle{\rm Cm}(\Delta(k)),\,{\bf
b}\rangle=\frac{n!}{(B_{(k)})^2}\Big(B_{(k)}\int_{F_{(k)j}}
\langle\mu_{(k)},\,{\bf b} \rangle  \, dX^j
-A_{(k)}\int_{F_{(k)j}} dX^j \Big).$$
 From (\ref{Ibr}) and (\ref{N(F)}), it follows
\begin{equation}\label{util}
\sum_{j=1}^m\frac{\partial}{\partial k_j}\,\langle{\rm
Cm}(\Delta(k)),\,{\bf b}\rangle=\frac{-1}{B_{(k)}}I(k;\,{\bf b}).
 \end{equation}
Thus, we have proved the following proposition:
 \begin{Prop}\label{iff}
$I(k;\,{\bf b})=0$ for all $k\in{\mathcal C}_{\Delta}\,$ iff
$\,\sum_{j=1}^m\frac{\partial}{\partial k_j}\,\langle{\rm
Cm}(\Delta(k)),\,{\bf b}\rangle=0$, for all $k\in{\mathcal
C}_{\Delta}$.
\end{Prop}

\smallskip

Next, we will parametrize the quotient ${\mathcal C}_{\Delta}/G$
(of classes of polytopes in ${\mathcal C}_{\Delta}$ module
translation) defined in the Introduction.

 After a possible renumbering, we may assume
that the intersection of facets $F_1,\dots,F_n$ is a vertex of
$\Delta$. Thus, the conormals ${\bf n}_1,\dots,{\bf n}_n$ are
linearly independent in ${\mathfrak t}$. So, given $k\in{\mathcal
C}_{\Delta}$, there is a unique $v\in{\mathfrak t}^*$, such that,
 \begin{equation}\label{defv}
 \langle v,\,{\bf n}_i\rangle+k_i=0,\;\;
i=1,\dots,n.
 \end{equation}
(Expressing the ${\bf n}_i$ in terms of a basis of ${\mathfrak t}$
and $v$ in the dual basis, (\ref{defv}) is a compatible and
determined system of linear equations for the coordinates of $v$.)
Moreover $v=v(k)$ depends {\em linearly} of the $k_i$; that is,
$\langle v(k),\,{\bf c}\rangle$ is a linear function of
$k_1,\dots,k_n$, for all ${\bf c}\in{\mathfrak t}$.

If $m-n=2$, we write
$$\lambda=k_{n+1}+\langle v(k),\,{\bf n}_{n+1}\rangle,\;\;\;\;
\tau=k_{m}+\langle v(k),\,{\bf n}_{m}\rangle,$$
  where $v(k)$ the element in ${\mathfrak t}^*$ defined by
  (\ref{defv}). From the linearity of $v$ with respect to the
  $k_i$, it follows that $\lambda$ and $\tau$ are {\em linear
  combinations} of $k_1,\dots,k_m$.

  The polytope in ${\mathcal C}_{\Delta}$ defined by
  $(k'_1=0,\dots,k'_n=0,\lambda,\tau)$ will be denoted by
  $\Delta_0(\lambda,\tau)$. It is the result of the
  translation of $\Delta(k)$ by the vector $v(k)$; i. e.,
 \begin{equation}\label{Cm(Delta0)}
  \Delta_0(\lambda,\tau)=\Delta(k)+v(k).
   \end{equation}

Let ${\bf b}$ an element in the integer lattice of ${\mathfrak
t}$, we define the function $g$ by
$$g(\lambda,\tau):=\langle{\rm
Cm}(\Delta_0(\lambda,\tau)),\,{\bf b}\rangle.$$
  The function $g$
is defined on the pairs $(\lambda,\tau)$ such that
 $(0,\dots,0,\lambda,\tau)\in{\mathcal C}_{\Delta}$.
 By Lemma \ref{LemCm}, it follows
$$g(s\lambda,s\tau)=sg(\lambda,\tau),$$
for any real number $s$ such that $(s\lambda,s\tau)$ belongs to
the domain of $g$. This property implies that
 \begin{equation}\label{glambda}
g=\lambda\frac{\partial g}{\partial \lambda}+\tau\frac{\partial
g}{\partial \tau}.
 \end{equation}

\begin{Thm}\label{vanishThm}
If $I(k; \,{\bf b})=0,$ for all $k\in{\mathcal C}_{\Delta}$ and
$r=2$, then
 $\langle {\rm Cm}(\Delta(k),\,{\bf
b}\rangle=\sum_j R_jk_j,$ with $R_j$
  constant (that is, $(\Delta,\,{\bf b})$ is a mass linear pair) and
 $\sum_jR_j=0$.
 \end{Thm}

{\it Proof.} We set $f(k_1,\dots,k_m):=\langle {\rm
Cm}(\Delta(k),\,{\bf b}\rangle$.
 It follows from (\ref{Cm(Delta0)}) that
 \begin{equation}\label{f=g-sum}
  f(k)=g(\lambda,\tau)- \langle v(k),\,{\bf b}\rangle.
  \end{equation}

By the hypothesis and Proposition \ref{iff},
 \begin{equation}\label{sumpartial=0}
 \sum_{j=1}^m\frac{\partial f}{\partial k_j}=0.
 \end{equation}

Since
 $$\sum_{j=1}^m\frac{\partial f}{\partial k_j}=
 \frac{\partial g}{\partial\lambda}\sum_{j=1}^m\frac{\partial \lambda}{\partial
 k_j}+
 \frac{\partial g}{\partial\tau}\sum_{j=1}^m\frac{\partial
\tau}{\partial k_j}-\big\langle \frac{\partial v}{\partial
k_j},\,{\bf b} \big\rangle,$$
 from (\ref{sumpartial=0}) we
deduce
 \begin{equation}\label{DiffEqug}
  p\frac{\partial
g}{\partial\lambda}+q\frac{\partial g}{\partial\tau}-t=0,
 \end{equation}
 where
$p,q,t$ stand for the following constants
$$p=\sum_{j=1}^m\frac{\partial \lambda}{\partial
 k_j},\;\;\;q=\sum_{j=1}^m\frac{\partial
\tau}{\partial k_j},\;\;\;t=\big\langle \frac{\partial v}{\partial
k_j},\,{\bf b} \big\rangle.$$

Since $q\lambda-p\tau$ and $t\tau-qg$ are first integrals of
(\ref{DiffEqug}), the general solution of this equation is
 \begin{equation}\label{g=Phi}
g(\lambda,\tau)=(t/q)\tau+\Phi(q\lambda-p\tau),
 \end{equation}
 where $\Phi$ is a derivable function of one variable.

It follows from (\ref{glambda}) and (\ref{g=Phi}) that
 \begin{equation}\label{Phi=uPhi}
\Phi(u)=u\Phi'(u).
 \end{equation}
  Thus, $\Phi(u)=\alpha u$, with $\alpha$
constant. We have for $f$
$$f(k)=(b/q)\tau+\alpha(q\lambda-p\tau)-\langle v(k),\,{\bf b}\rangle.$$
 In other words, $f$ is a linear function of the $k_j$; i.e.,
 $f(k)=\sum_jR_jk_j$, with $R_j$ constant.
 From (\ref{sumpartial=0}), it follows $\sum_jR_j=0$.
 \qed

 \smallskip

{\it Remark.} The proof of Theorem \ref{vanishThm} can be adapted
to the simpler case when $r=1$. In this case, the function
$g(\lambda)=\langle{\rm Cm}(\Delta_0(\lambda),\,{\bf b}\rangle$,
satisfies $p\frac{{\rm d}g}{{\rm d}\lambda}-t=0$ and
$g(s\lambda)=sg(\lambda)$. So, $g(\lambda)=(t/p)\lambda$ and
$f(k)=(t/p)\lambda+\langle v(k),\,{\bf b}\rangle$ is a linear map
of the variables $k_j$.

 On the other hand,  the proof of this theorem does not
admit an adaptation  to the case   $r>2$. In fact,   the
corresponding function $\Phi$ would be  a function of $r-1$
variables $\Phi(u_1,\dots,u_{r-1})$. The equation which
corresponds to (\ref{Phi=uPhi}) in this case  would be
$$\Phi=\sum_{i=1}^{r-1}u_i\frac{\partial \Phi}{\partial u_i}.$$
 But this condition does not implies the
linearity of $\Phi$.

\medskip

When $(\Delta,\,{\bf b})$ is a mass linear pair as in
(\ref{Rjkj}), by (\ref{util})
 \begin{equation}\label{IBsumR}
I(k;\,{\bf b})=-B_{(k)}\sum_jR_j,
 \end{equation}
 for all $k\in{\mathcal
C}_{\Delta}.$ From  (\ref{IBsumR}), we deduce  the following
proposition:
\begin{Prop}\label{Corofinal}
Let $(\Delta,\,{\bf b})$ be a mass linear pair. $I(k;\,{\bf b})=0$
for all $k\in{\mathcal C}_{\Delta}$ iff $\sum_jR_j=0$.
 \end{Prop}

\medskip

 {\bf Proof of Theorem \ref{FinalThm}.}
 It is a direct consequence of Proposition \ref{Corofinal},  Theorem
 \ref{vanishThm} and the Remark above.
  \qed

\medskip

We will deduce a sufficient condition for a Delzant polytope
$\Delta$ to  admit mass linear functions. We write
$$\Dot{\rm Cm}(\Delta(k)):=\frac{{\rm d}}{{\rm
d}\,\epsilon}\bigg|_{\epsilon=0}{\rm
Cm}(\Delta(k+\check\epsilon)),$$
 with $\check\epsilon=(\epsilon,\dots,\epsilon).$

\begin{Prop}\label{PropSufCon}
If all  points $\Dot{\rm Cm}(\Delta(k))$, for $k\in{\mathcal
C}_{\Delta}$, belong to a hyperplane of $({\mathbb R}^n)^*$ with a
conormal vector in ${\mathbb Z}^n$ and $r\leq 2$, then $\Delta$
admits a mass linear function.
\end{Prop}

{\it Proof.} Let ${\bf b}\in{\mathbb Z}^n$ be a conormal vector to
the hyperplane, then
$$0=\langle \Dot{\rm Cm}(\Delta(k)),\,{\bf b}
\rangle=\big\langle\sum_j\frac{\partial}{\partial k_j}\,{\rm
Cm}(\Delta(k)) ,\,{\bf b}\big\rangle.$$
 By (\ref{util}), $I(k,\,{\bf b})=0$; Theorem
 \ref{vanishThm} applies and $(\Delta,\,{\bf b})$ is a mass linear pair.
  \qed

\smallskip

\begin{Prop}\label{PropNecCond}
 Let $\Delta$ be a Delzant polytope, such that $k=0$ belongs to the
 closure of ${\mathcal C}_{\Delta}$. If $r\leq 2$, a necessary condition for
 the vanishing of $I(k;\,{\bf b})$ on ${\mathcal C}_{\Delta}$ is
  $$\big\langle\frac{{\rm d}}{{\rm d}\epsilon}\bigg|_{\epsilon=0}{\rm
 Cm}(\Delta(\check\epsilon)),\,{\bf b}\big\rangle=0.$$
 \end{Prop}
  {\it Proof.} If $I(k;\,{\bf b})$ vanishes on ${\mathcal C}_{\Delta}$, then
    $(\Delta,\,{\bf b})$ is a linear
  pair, by Theorem \ref{FinalThm}.  Thus,
  $\langle {\rm Cm}(\Delta(k)),\,{\bf b}\rangle=\sum_jR_jk_j+C$,
  on ${\mathcal C}_{\Delta}$.   So, given $k\in{\mathcal C}_{\Delta}$ and $\epsilon$ small
  enough
  $$\langle {\rm Cm}(\Delta(k+\check \epsilon)),\,{\bf b}\rangle=\sum_jR_jk_j+\epsilon\sum_jR_j +C.$$
  By  Theorem \ref{FinalThm},  $\sum_j R_j=0$. Thus, for any
  $k\in{\mathcal C}_{\Delta}$,
  $$\frac{{\rm d}}{{\rm d}\epsilon}\bigg|_{\epsilon=0} \langle {\rm Cm}(\Delta(k+\check \epsilon)),\,{\bf b}\rangle =0.$$
Taking the limit as $k\to 0$,
$$0=\lim_{k\to 0}\frac{{\rm d}}{{\rm d}\epsilon}\bigg|_{\epsilon=0} \langle {\rm Cm}(\Delta(k+\check \epsilon)),\,{\bf
b}\rangle=\big\langle\frac{{\rm d}}{{\rm
d}\epsilon}\bigg|_{\epsilon=0}{\rm
Cm}(\Delta(\check\epsilon)),\,{\bf b}\big\rangle.$$
 \qed

\smallskip

 Next, we will describe a geometric interpretation of the number $\sum_jR_j$.
  Given  an arbitrary Delzant polytope $\Delta$.
 If $a$ is a vector of ${\mathfrak t}^*$, then
 \begin{equation}\label{CmTrnaslation}
  {\rm Cm}(\Delta(k'))={\rm Cm}(\Delta(k))+a,
  \end{equation}
if $k'_j=k_j+\langle a,\,  {\bf n}_j\rangle.$

We will denote by $d$  the element of ${\mathfrak t}^*$ defined by
the following relation
\begin{equation}\label{kj+1}
 {\rm Cm}(\Delta(\tilde k))={\rm Cm}(\Delta(k))+d,
 \end{equation}
 with $\tilde k_j=k_j+1$ for all $j$.

From (\ref{CmTrnaslation}) and (\ref{kj+1}),  we have
 $${\rm Cm}(\Delta(k_j+\langle d,\, {\bf n}_j\rangle))=
{\rm Cm}(\Delta(k_j))+ d=
 {\rm Cm}(\Delta(\tilde k_j=k_j+1)).$$

 Now, we assume that  $(\Delta,\,{\bf b})$ is a mass linear pair.
 % as in (\ref{util}).
  From (\ref{Rjkj}),  it follows
 $$\langle {\rm Cm}(\Delta(k_j+\langle d,\, {\bf n}_j\rangle)),\,{\bf
 b}\rangle=\sum R_jk_j+\sum R_j\langle d,\, {\bf n}_j\rangle+C.$$
  $$\langle {\rm Cm}(\Delta(k_j))+ d,\,{\bf
 b}\rangle=\sum R_jk_j+ \langle d,\,{\bf b}\rangle+C,\;\;
 \langle {\rm Cm}(\Delta(\tilde k)),\,{\bf
 b}\rangle=\sum R_jk_j+ \sum R_j+C.$$
These formulas allow us to state the following proposition, that
gives an interpretation of the sum $\sum_jR_j$ in terms of the
variation of ${\rm Cm}(\Delta(k))$ with the $k_j$.
\begin{Prop}\label{PropsumR}
 Let $(\Delta,\,{\bf b})$ be a mass linear pair as in (\ref{Rjkj}).
 Then,
 $$\sum_j R_j\langle d,\, {\bf n}_j\rangle = \langle d,\, {\bf
b}\rangle = \sum_j R_j,$$
  $d$ being the element of ${\mathfrak t}^*$ defined by
  (\ref{kj+1}).
\end{Prop}

%%%%%%%%%%%%%%%%%%%%%%%%%%%%%%%%%%%%%%%%%%%%%%%%%%%%%%%%%%%%%%%%%%%%%%%%%%%%%%%%%%%%%%%%%%%%%%%%%%%%%%%%%%%%%%%%%%
%%%%%%%%%%%%%%%%%%%%%%%%%%%%%%%%%%%%%%%%%%%%%%%%%%%%%%%%%%%%%%%%%%%%%%%%%%%%%%%%%%%%%%%%%%%%%%%%%%%%%%%%%%%%%%%%%%%%%%%

\section{Examples}\label{Aplications}

In this Section, we will deduce the particular form which adopts
Theorem
 \ref{CorThm1}, when $\Delta$ is a  polytope of the
 types (i)-(iii) mentioned in the Introduction. For each case,
 we will determine the center of mass of the corresponding polytope $\Delta(k)$ and the condition
 for $(\Delta,\,{\bf b})$ to be a mass linear pair. We will
 dedicate a subsection to each type.

\subsection{Hirzebruch surfaces}

Given $r\in{\Bbb Z}_{>0}$ and $\tau,\lambda\in{\Bbb R}_{>0}$ with
$\sigma:=\tau-r\lambda>0$, in \cite{Vipp} we considered the
Hirzebruch surface $N$ determined by these numbers. $N$ is the
quotient
$$\big\{z\in {\Bbb C}^4\,\,:\,\, |z_1|^2+r|z_2|^2+|z_4|^2=\tau/\pi,\;
|z_2|^2+|z_3|^2=\lambda/\pi\big\}/{\Bbb T},$$ where the
equivalence defined by ${\Bbb T}=(S^1)^2$ is given by
$$(a,b)\cdot(z_1,z_2,z_3,z_4)=(az_1, a^rbz_2,bz_3,az_4),$$
for $(a,b)\in (S^1)^2$.

The manifold $N$ equipped with the following $(U(1))^2$ action
$$(\epsilon_1,\,\epsilon_2)[z_j]=[\epsilon_1 z_1,\,\epsilon_2
z_2,\,z_3,\,z_4],$$ is a toric manifold. The corresponding moment
 polytope $\Delta$ is the trapezium in ${\mathbb R}^2$ with
vertices
 \begin{equation}\label{vertices}
P_1=(0,0),\; \; P_2=(0,\lambda),\;\; P_3=(\tau,0),\;\;
P_4=(\sigma,\lambda).
 \end{equation}
 That is, $N$ is the toric manifold $M_{\Delta}$ determined by the trapezium $\Delta$.

As the conormals to the facets of $\Delta$ are the vectors ${\bf
n}_1=(-1,0)$, ${\bf n}_2=(0,-1)$ ${\bf n}_3=(0,1)$ and ${\bf
n}_4=(1,r)$, the facets of a generic polytope $\Delta(k)$ in
${\mathcal C}_{\Delta}$ are on the straights
$$-x=k_1,\;\;-y=k_2,\;\;y=k_3,\;\;x+ry=k_4.$$
The vertices of $\Delta(k)$ are the points
$$(-k_1,\,-k_2),\; \;(-k_1,\,k_3), \;\;(k_4-rk_3,\,k_3),
\;\;(k_4+rk_2,\,-k_2).$$
 Thus, the translation in the plane $x,y$ defined by $(-k_1,\,-k_2)$ transforms the trapezium
 determined by the vertices (\ref{vertices}) in $\Delta(k)$, if
 \begin{equation}\label{tau=k}
  \tau=k_4+rk_2+k_1,\;\;\lambda=k_3+k_2.
  \end{equation}
   So,
 \begin{equation}\label{Cm+k}
  {\rm Cm}(\Delta(k))={\rm Cm}(\Delta)+(-k_1,\,-k_2).
  \end{equation}
 Moreover, the mass center of $\Delta$ is
\begin{equation}\label{Cmblowup}{\rm
Cm}(\Delta)=\Big(\frac{3\tau^2-3r\tau\lambda+r^2\lambda^2}{3(2\tau-r\lambda)},\;
\frac{3\lambda\tau-2r\lambda^2}{3(2\tau-r\lambda)} \Big).
\end{equation}

  The chamber ${\mathcal C}_{\Delta}$ consists of the points $(k_1,\dots, k_4)$ such that
  $\tau-r\lambda>0$, with $\tau$ and $\lambda$ given by
  (\ref{tau=k}). So,
   the point $k=0$ belongs to the closure of ${\mathcal
   C}_{\Delta}$.
   From (\ref{Cm+k}), together with (\ref{tau=k}) and
   (\ref{Cmblowup}), it follows
   \begin{equation}\label{Cm(check)}
   {\rm Cm}(\Delta(\check\epsilon))=\big( \frac{r^2\epsilon}{12},\,\frac{-r\epsilon}{6}
   \big),
    \end{equation}
     where $\check\epsilon=(\epsilon,\epsilon,\epsilon,\epsilon)$.
   By Proposition \ref{PropNecCond}, if $I(k;\,{\bf b})$ with ${\bf b}=(b_1,\,b_2)\in{\mathbb Z}^2$ vanishes
   on the chamber ${\mathcal C}_{\Delta}$,  then $rb_1-2b_2=0$.

\smallskip

 On the other hand, from (\ref{Cmblowup}) and (\ref{Cm+k}), it follows
\begin{equation}\label{CmShort}
 \langle {\rm Cm}(\Delta(k)),\,{\bf b}\rangle=\frac{(
3\tau^2-3r\tau\lambda+r^2\lambda^2  )b_1+(3\lambda\tau-2r\lambda^2
)b_2 }{3(2\tau-r\lambda) }-k_1b_1-k_2b_2.
 \end{equation}
 By (\ref{tau=k}),  expression (\ref{CmShort})  is linear in the $k_i$ iff
 $$\frac{(
3\tau^2-3r\tau\lambda+r^2\lambda^2  )b_1+(3\lambda\tau-2r\lambda^2
)b_2 }{3(2\tau-r\lambda) }$$ is linear in $\tau,\lambda$. That is,
iff there exist constants $A,B$ such that for al $\tau,\lambda$
$$(
3\tau^2-3r\tau\lambda+r^2\lambda^2  )b_1+(3\lambda\tau-2r\lambda^2
)b_2=3(2\tau-r\lambda)(A\tau+B\lambda).$$
 From this relation, it
follows the above condition   $rb_1=2b_2$. In this case
(\ref{CmShort}) reduces to
\begin{equation}\label{linearHirz}
\langle{\rm Cm}(\Delta(k),\,{\bf
b}\rangle=\frac{-b_1}{2}k_1+\frac{b_1}{2}k_4.
 \end{equation}

Comparing (\ref{Rjkj}) with   (\ref{linearHirz}), we obtain,
$R_1=-R_4=\frac{-b_1}{2}$, $\,R_2=R_3=0$;  so,  $\sum_jR_j=0$.
 That is, the condition $\sum_jR_j=0$ holds for all the mass pairs
 $(\Delta,\,{\bf b})$ when $\Delta$ is the polytope associated to
 a Hirzebruch surface.
 Hence,
 \begin{Prop}\label{PropHirz}
$(\Delta,\,{\bf b})$ is a mass linear pair iff $rb_1=2b_2$.
Moreover, in  this case $\sum_jR_j=0$.
\end{Prop}

By Theorem \ref{CorThm1}, we have

 \begin{Cor}\label{CoroHirz}
If $rb_1\ne 2b_2$, then $\psi_{\bf b}$ generates an infinite
cyclic subgroup in $\pi_1({\rm
Ham}(M_{\Delta},\,\omega_{\Delta}))$.
\end{Cor}

\smallskip

{\it Remark.}

  We denote by $\phi_t$ the following isotopy of $M_{\Delta}$
$$\phi_t[z]=[e^{2\pi it}z_1,z_2,z_3,z_4].$$
$\phi$ is a loop in the Hamiltonian group of $M_{\Delta}$. By
$\phi'$ we denote the Hamiltonian loop
$$\phi'_t[z]=[z_1,e^{2\pi it}z_2,z_3,z_4].$$
In Theorem 8 of \cite{V} we proved that
$I({\phi'})=(-2/r)I(\phi)$. If  ${\bf b}=(b_1,\,b_2)\in{\mathbb
Z}^2$, then
$$I({\psi_{\bf b}})=b_1I({\phi})+b_2I({\phi'})=(b_1-(2/r)b_2)I({\phi}).$$
 That is, $I(\psi_{\bf b})=0$ iff $rb_1=2b_2$, which is in agreement with
 Proposition \ref{PropHirz} and Theorem \ref{FinalThm}.

%%%%%%%%%%%%%%%%%%%%%%%%%%%%%%%%%%%%%%%%%%%%%%%%%%%%%%%%%%%%%%%%%%%%%%%%%%%%%%%%%%%%%%%%%%%%%%%%%%%%%%%%

 \medskip

 \subsection{$\Delta _p$ bundle over $\Delta_1$}

Given the integer $p>1$, as McDuff and Tolman in
%following
\cite{M-T3}, we consider the following vectors in ${\mathbb
R}^{p+1}$
  \begin{equation}\label{definition-ni}
   {\bf n}_i=-e_i,\; i=1,\dots,p,\;\;{\bf n}_{p+1}=\sum_{i=1}^p
  e_i,\;\; {\bf n}_{p+2}=-e_{p+1},\;\; {\bf
 n}_{p+3}=e_{p+1}-\sum_{i=1}^pa_ie_i,
  \end{equation}
  where $e_1,\dots, e_{p+1}$
is the standard basis of ${\mathbb R}^{p+1}$ and $a_i\in{\mathbb
Z}.$ We write
$${\bf a}:=(a_1,\dots,a_p)\in{\mathbb Z}^{p},\;\;A:=\sum_{i=1}^p
a_i,\;\;{\bf a}\cdot {\bf a}=\sum_{i=1}^p a_{i}^2.$$

 Let  $\lambda,\,\tau$  be real positive numbers with $\lambda+a_i>0,$
for $i=1,\dots,p$. In this subsection we will consider the
polytope $\Delta$ in $({\mathbb R}^{p+1})^*$ defined by the above
conormals ${\bf n}_j$ and the following $k_j$
 \begin{equation}\label{constki}
 k_1=\dots=k_p=k_{p+2}=0,\; k_{p+1}=\tau,\;
k_{p+3}=\lambda.
 \end{equation}

This polytope will be also denote by $\Delta_0(\lambda,\tau)$. It
is a $\Delta_p$ bundle on $\Delta_1$ (see \cite{M-T3}). When
$p=2$, $\Delta=\Delta_0(\lambda,\tau)$ is the prism whose base is
the triangle of vertices $(0,0,0),$ $(\tau,0,0)$ and $(0,\tau,0)$
and whose ceiling is the triangle determined by $(0,0,\lambda),$
$(\tau,0,\lambda+a_1\tau)$ and $(0,\tau,\lambda+a_2\tau)$.

%(see Figure 1).

%%%%%%%%%%%%%%%%%%%%%%%%%%%%%%%%%%%%%%%%%%%%%%%%%%%%%%%%%%%%%%%%%%%%%%%%%%%%%%%%%%%%%%%%%%%%%%%%%%%%%%%%%%%%%%%%
%%%%%%%%%%%%%%%%%%%%%%%%%%%%%%%%%%%%%%%%%%%%%%%%%%%%%%%%%%%%%%%%%%%%%%%%%%%%%%%%%%%%%%%%%%%%%%%%%%%%%%%%%%
%\begin{figure}[htbp]
%\begin{center}
%\epsfig{file=Polytope1.eps, height=4cm}
% \caption[Figure 1]{\small $\Delta_2$ bundle over $\Delta_1$}
%\end{center}
%\end{figure}
%%%%%%%%%%%%%%%%%%%%%%%%%%%%%%%%%%%%%%%%%%%%%%%%%%%%%%%%%%%%%%%%%%%%%%%%%%%%%%%%%%%%%%%%%%%%%%%%%%%%%%%%
%%%%%%%%%%%%%%%%%%%%%%%%%%%%%%%%%%%%%%%%%%%%%%%%%%%%%%%%%%%%%%%%%%%%%%%%%%%%%%%%%%%%%%%%%%%%%%%%%%%%%%%%%%%%%%%%%%%%%%%%%%%%%%555555

We assume that the above polytope $\Delta$ is a Delzant polytope.
The manifold (\ref{DefiM}) is in this case
$$M_{\Delta}=
 \big\{ z\in{\mathbb C}^{p+3}\,:\, \sum_{i=1}^{p+1}|z_i|^2=\tau/\pi,\;
\;-\sum_{j=1}^p a_j|z_j|^2+|z_{p+2}|^2+|z_{p+3}|^2=\lambda/\pi
    \big\}/ \simeq,$$ where $(z_j)\simeq (z'_j)$ iff there are
$\alpha,\beta\in U(1)$ such that
$$z'_j=\alpha\beta^{-a_j}z_j,\;j=1,\dots,p\,; \; \;\; z'_{p+1}=\alpha z_{p+1};\;\;\;z'_k=\beta z_k,\; k=p+2,\,p+3.$$
Thus, $M_{\Delta}$  is the total space of the fibre bundle
${\mathbb P}(L_1\oplus\dots \oplus L_p\oplus{\mathbb C})\to
{\mathbb C}P^1$, where $L_j$ is the holomorphic line bundle over
${\mathbb C}P^1$ with Chern number $a_j$.

The symplectic form (\ref{omeg}) is
$$\omega_{\Delta}=(1/2)\big(\sigma_1+\dots+\sigma_p+\sigma_{p+2}),$$
where $\sigma_k= d\rho_k^2\wedge d\varphi_k.$

 And the moment map
\begin{equation}\label{momentb}
\mu_{\Delta}([z])=(x_1,\dots, x_{p},\,x_{p+2}),
\end{equation}
where $x_i:=\pi\rho_i^2.$

\smallskip

\begin{Prop}\label{CmDeltabundle}
The coordinates $\bar x_j$ of ${\rm Cm}(\Delta_0(\lambda,\tau))$
are given by:
 \begin{equation}\label{CmbMT}
\bar{x}_k=\frac{\tau}{p+2}\,  \frac{\lambda(p+2)+\tau\big(A+ a_k
\big)}
 {\lambda(p+1)+\tau A},\;\; \hbox{for}\;\; k=1,\dots, p.
\end{equation}
 \begin{equation}\label{CmbMT1}
 \bar{x}_{p+2}= \frac{1}{2}\,\frac{(p+1)(p+2)\lambda^2+2(p+2)A\lambda\tau+({\bf
a}\cdot {\bf a}+A^2)\tau^2}{(p+2)\big((p+1)\lambda+A \tau \big)}.
\end{equation}
\end{Prop}

{\it Proof.} Since the points $[z]\in M_{\Delta}$ satisfy
$|z_{p+2}|^2\leq\lambda/\pi+\sum_{j=1}^pa_j|z_j|^2,$ by
(\ref{Intpolytope1}) and Lemma \ref{Lemmaint} we have
 \begin{equation}\label{intMomega}
  \int_{M_{\Delta}}(\omega_{\Delta})^{p+1}=(p+1)!\int_{S_p(\tau)}\big(\lambda+\sum_{j=1}^p
a_jx_j\big)=(p+1)!\Big(\frac{\lambda\tau^p}{p\,!}+\frac{\tau^{p+1}A}{(p+1)!}\Big).
 \end{equation}
Similarly, for $k=1,\dots, p$
\begin{equation}\label{itnxkOmega}
 \int_{M_{\Delta}}x_k\,(\omega_{\Delta})^{p+1}=(p+1)!\Big(\frac{\lambda\tau^{p+1}}{(p+1)!}+
 \frac{\tau^{p+2}}{(p+2)!}\sum_{j\ne k}a_j+  \frac{2\tau^{p+2}a_k}{(p+2)!}
 \Big).
  \end{equation}
     The $k$-th coordinate of $\text{Cm}(\Delta)$, $\bar{x}_{k}$,
  is the quotient of (\ref{itnxkOmega})  by (\ref{intMomega}); that
  is,
 $$\bar{x}_k=\frac{\tau}{p+2}\,  \frac{\lambda(p+2)+\tau\big(A+ a_k   \big)}
 {\lambda(p+1)+\tau A}.$$

For the $p+2$-coordinate of ${\rm Cm}(\Delta)$, we need
 to calculate $\int_M x_{p+2}(\omega_{\Delta})^{p+1}.$
By Lemma \ref{Lemmaint}
 \begin{align} \label{alignCM}
  \frac{1}{(p+1)!}\int_M
x_{p+2}(\omega_{\Delta})^{p+1}=&\frac{1}{2}\int_{S_p(\tau)}\big(\lambda+\sum_{j=1}^pa_jx_j
\big)^2 \\ \notag
 =&\frac{1}{2}\Big(\frac{\lambda^2\tau^p}{p\,!}+
\frac{2A\lambda\tau^{p+1}}{(p+1)!}+
 \frac{({\bf a}\cdot{\bf a}+A^2)\tau^{p+2}}{(p+2)!} \Big).
 \end{align}
Formula (\ref{CmbMT1}) is a consequence of (\ref{intMomega})
together with (\ref{alignCM}).
 \qed

\smallskip

The translation in $({\mathbb R}^{p+1})^*$ defined by the vector
$(-k_1,\dots,-k_p,-k_{p+2})$ transforms the hyperplanes $\langle
x,\,{\bf n}_{p+3}\rangle=\lambda$ and $\langle x,\,{\bf
n}_{p+1}\rangle=\tau$ in
 \begin{equation}\label{translw}
 \langle x,\,{\bf n}_{p+3}\rangle=\lambda-k_{p+2}+\sum_{j=1}^pa_jk_j,\;\;\;\; \langle x,\,{\bf
n}_{p+1}\rangle=\tau-\sum_{j=1}^p k_j,
 \end{equation}
respectively.

Let $\Delta(k)$ be a polytope with $k=(k_1,\dots,k_{p+3})$ generic
in the chamber ${\mathcal C}_{\Delta}$. From (\ref{translw}),  it
follows that $\Delta(k)$ is the image of the polytope
$\Delta_0(\lambda,\tau)$
 by the
translation determined by  $(-k_1,\dots,-k_p,-k_{p+2})$, whenever
\begin{equation}\label{k-lamb-tau}
k_{p+2}-\sum_{j=1}^pa_jk_j+k_{p+3}=\lambda,\;\;\;\,\sum_{j=1}^pk_j+k_{p+1}=\tau.
 \end{equation}
 In this case,
 \begin{equation}\label{Cm(D(k))}
 {\rm Cm}(\Delta(k))={\rm
 Cm}(\Delta_0(\lambda,\tau))-(k_1,\dots,k_p,k_{p+2}).
 \end{equation}

 According to (\ref{k-lamb-tau}), the coordinates of the  mass center ${\rm Cm}(\Delta(\check\epsilon))$, with
 $\check\epsilon=(\epsilon,\dots,\epsilon)$, can
  be  obtained
 substituting in (\ref{CmbMT}) and in (\ref{CmbMT1})  $\lambda$ by
 $$\epsilon-\sum_{j=1}^p a_j\,\epsilon+\epsilon=(2-A)\epsilon$$
 and $\tau$ by $(p+1)\epsilon$, and finally   take into account (\ref{Cm(D(k))}).  These operations give
  \begin{align}\notag
 &\bar{x}_j(\Delta(\check \epsilon))=\frac{\epsilon}{2(p+2)}\big(
(p+1)a_j-A \big),
 \;\;\; j=1,\dots,p \\ \notag
 &\bar{x}_{p+2}(\Delta(\check
 \epsilon))=\frac{\epsilon}{4(p+2)}\big(-A^2+(p+1)({\bf a}\cdot{\bf
 a})\big).
  \end{align}

 Given ${\bf b}=(b_1,\dots,b_p,b)\equiv({\bf\Hat b},\bf{\Dot b})$,
 with ${\bf\Hat b}=(b_1,\dots,b_p,0)$ and ${\bf\Dot
 b}=(0,\dots,0,b)$.
$$ \big\langle\frac{{\rm d}}{{\rm d}\epsilon}\bigg|_{\epsilon=0}{\rm
 Cm}(\Delta(\check\epsilon)),\,{\bf b}\big\rangle=
 \frac{1}{4(p+2)}\Big((p+1)\big( 2\,{\bf a}\cdot{\bf\Hat b}  -b\,{\bf a}\cdot{\bf a}\big)-A(2B+bA)
 \Big),$$
 where ${\bf a}\cdot{\bf\Hat b}=\sum_{j=1}^p a_jb_j$ and $B=\sum_{j=1}^p b_j.$

 By Proposition \ref{PropNecCond}, we have:

 \begin{Thm}\label{pi1fifration}
 Let $\Delta$ be the $\Delta_p$ bundle over $\Delta_1$ defined by
 (\ref{definition-ni}) and (\ref{constki}). Given ${\bf b}=(\Hat{\bf b},\Dot{\bf
 b})\in{\mathbb Z}^{p+1}$, if
 $$(p+1)\big( 2{\bf a}\cdot{\bf\Hat b} -b\,{\bf a}\cdot{\bf
 a}\big)-A(2B+bA)\ne 0,$$
 then $\psi_{\bf b}$ defines an infinite cyclic subgroup in the
 fundamental group $\pi_1({\rm Ham}(M_{\Delta},\,\omega_{\Delta})).$
 \end{Thm}

It is straightforward to check that
 \begin{equation}\label{straightforward}
  (p+1)\big( 2{\bf a}\cdot{\bf \Hat b} -b\,{\bf a}\cdot{\bf
 a}\big)-A(2B+bA)=0
 \end{equation}
  is also a sufficient condition for $(\Delta,{\bf b})$ to be a mass linear
 pair.
% By Theorem \ref{FinalThm} this condition holds iff $I(\Delta,{\bf b})$ vanishes.
 %   In
% \cite{Vipp} Theorem 11, this result
% has been proved by direct calculation.

\smallskip

Since $$\langle {\rm Cm} (\Delta_0(\lambda,\tau)),\,{\bf
b}\rangle=
 \langle {\rm Cm} (\Delta_0(\lambda,\tau)),\,{\bf\Hat b}\rangle +
 \langle {\rm Cm} (\Delta_0(\lambda,\tau)),\,{\bf\Dot b}\rangle,$$
 if (\ref{straightforward}) holds, using (\ref{CmbMT}) and
 (\ref{CmbMT1}), one obtains
 $$\langle
{\rm Cm} (\Delta_0(\lambda,\tau)),\,{\bf
b}\rangle=\frac{b\lambda}{2}+ \Big( \frac{b}{2}\,\frac{ ({\bf
a}\cdot{\bf a}+A^2  )} {(p+2)A}+ \frac{ ({\bf a}\cdot{\bf \Hat
b}+AB  )} {(p+2)A}\Big)\tau.$$
 By (\ref{Cm(D(k))}), for $k\in{\mathcal C}_{\Delta}$,
$$ \langle
{\rm Cm} (\Delta(k)),\,{\bf b}\rangle =\langle {\rm Cm}
(\Delta_0(\lambda,\tau)),\,{\bf
b}\rangle-\sum_{j=1}^pb_jk_j-bk_{p+2},$$
 with $\lambda$ and $\tau$ given by (\ref{k-lamb-tau}).

If ${\bf b}={\bf \Hat b}$, the condition (\ref{straightforward})
reduces to $(p+1){\bf a}\cdot{\bf \Hat b}=AB$ and
 $$\langle
{\rm Cm} (\Delta(k)),\,{\bf b}\rangle=\frac{ ({\bf a}\cdot{\bf
\Hat b}+AB  )} {(p+2)A}\Big( \sum_{j=1}^pk_j+k_{p+1}
\Big)-\sum_{j=1}^pb_jk_j.$$
 Hence, $\langle
{\rm Cm} (\Delta(k)),\,{\bf b}\rangle=\sum R_jk_j$, where
$$R_j=\frac{ ({\bf a}\cdot{\bf
\Hat b}+AB  )} {(p+2)A}-b_j,\;\;j=1,\dots,p;\;\;\;R_{p+1}=\frac{
({\bf a}\cdot{\bf \Hat b}+AB  )} {(p+2)A},\;\;
R_{p+2}=R_{p+3}=0.$$
 So,
 $$\sum_{j=1}^{p+3}R_j=\frac{(p+1){\bf a}\cdot{\bf \Hat
 b}-AB}{(p+2)A}=0.$$

A similar calculation for the case ${\bf b}={\bf \Dot b}$ shows
that the corresponding $\sum_jR_j$ vanishes. That is,
 \begin{Prop}\label{SumRjBundle}
Let $\Delta$ be a $\Delta_p$ bundle over $\Delta_1$. If
$(\Delta,{\bf b})$ is a mass linear pair, then $\sum_jR_j=0$.
 \end{Prop}

 \smallskip

For $p=2$, let ${\bf b}$ be the following linear combination of
the conormal vectors ${\bf b}=\gamma_1{\bf n}_1+\gamma_2{\bf
n}_2+\gamma_3{\bf n}_3$ with $\gamma_1+\gamma_2+\gamma_3=0$. By
(\ref{definition-ni}), ${\bf b}=(b_1,b_2,0)$ with
$b_1=\gamma_3-\gamma_1$, $b_2=\gamma_3-\gamma_2$. In this case
condition (\ref{straightforward}) reduces to
$$3(a_1b_1+a_2b_2)=(a_1+a_2)(b_1+b_2).$$
Or in terms of the $\gamma_i$
 \begin{equation}\label{p=2}
 a_1\gamma_1+a_2\gamma_2=0.
 \end{equation}
 This is a necessary and sufficient condition for $(\Delta,{\bf
 b})$ to be mass linear. This result is the statement of Lemma 4.8
 in \cite{M-T3}.

\medskip

%%%%%%%%%%%%%%%%%%%%%%%%%%%%%%%%%%%%%%%%%%%%%%%%%%%%%%%%%%%%%%%%%%%%%%%%%%%%%%%%%%%%%%%%%%%%%%%%%%%%%%%%%%%

\subsection{One point blow up of ${\mathbb C}P^n$.} \label{SubsectCP2}

 In this subsection $\Delta\equiv\Delta_0(\lambda,\tau)$ will be
\begin{equation}\label{Deltaonepoint}
\Delta=\Big\{(x_1,\dots,x_n)\in{\mathbb R}^n\,|\,
\sum_{i=1}^nx_i\leq\tau,\;0\leq x_i,\; x_n\leq\lambda\Big\},
\end{equation}
 where $\tau,\lambda\in{\mathbb R}_{>0}$ and
$\sigma:=\tau-\lambda>0$. That is, $\Delta$ is the polytope
obtained truncating the simplex $S_n(\tau)$, defined in Lemma
\ref{Lemmaint}, by a ``horizontal" hyperplane through the point
$(0,\dots,0,\lambda).$
 The manifold $M_{\Delta}$ associated with $\Delta$ is the one
point blow up of
 ${\mathbb C}P^n$.

 The mass center of the simplex
$S_n(\tau)$ is the point
 \begin{equation}\label{CmSn}{\rm
Cm}(S_n(\tau))=\frac{\tau}{n+1}w,
 \end{equation}
  with $w=(1,\dots,1)$.

As the volume   of $S_n(\tau)$ is $\tau^n/n!$, it follows from
(\ref{CmSn})
$$(\tau^n-\sigma^n)\,{\rm
Cm}(\Delta)=\tau^n\frac{\tau}{n+1}w-\sigma^n\big(\frac{\sigma}{n+1}w+\lambda\,
e_n \big).$$ That is,
 \begin{equation}\label{Cmtruncated}
  {\rm Cm}(\Delta)=
\frac{1}{\tau^n-\sigma^n}\Big(\big(\frac{\tau^{n+1}-\sigma^{n+1}}{n+1}
\big)w-\lambda\sigma^ne_n \Big).
\end{equation}

Given $k=(k_1,\dots,k_{n+2})\in{\mathcal C}_{\Delta}$, the facets
of $\Delta(k)$ are in the following hyperplanes:
\begin{equation}\label{hyperblow}
-x_j=k_j,\,j=1,\dots,n;\;\;\;\sum_{i=1}^px_k=k_{n+1};\;\;\;x_{n+1}=k_{n+2}.
 \end{equation}
As in the preceding subsections,
% for $k\in{\mathcal C}_{\Delta}$,
 \begin{equation}\label{Delta(k)Blow}
 \Delta(k)=\Delta_0(\lambda,\tau)-(k_1,\dots,k_n),
  \end{equation}
  provided
$\lambda=k_n+k_{n+2}$ and $\tau=\sum_{i=1}^{n+1}k_i.$

The pair $(\Delta,\, {\bf b}=(b_1,\dots, b_n))$ is mass linear iff
there exist $A,B,C\in{\mathbb R}$ such that
$$\sum_{j=1}^{n-1}b_j\frac{\tau^{n+1}-\sigma^{n+1}}{n+1}+b_n\Big(
\frac{\tau^{n+1}-\sigma^{n+1}}{n+1}-(\tau-\sigma)\sigma^n
 \Big)=(A\tau+B\sigma+C)(\tau^n-\sigma^n),$$
 for all $\tau,\sigma$ ``admissible". A simple calculation proves  the
 following proposition:
 \begin{Prop}\label{PropLinear}
The pair $(\Delta, {\bf b})$ is mass linear iff
$$b_n=\frac{1}{n}\sum_{j=1}^{n-1}b_j.$$
\end{Prop}

 From Theorem \ref{CorThm1} together with Proposition
\ref{PropLinear}, it follows the following theorem:

 \begin{Thm}\label{Propblowup}
 If ${\bf b}=(b_1,\dots,b_n)\in{\mathbb Z}^n$ and
$\sum_{j=1}^{n-1}b_j\ne nb_n$, then $\psi_{\bf b}$ generates an
infinite cyclic subgroup in
$\pi_1(\rm{Ham}(M_{\Delta},\,\omega_{\Delta}))$.
\end{Thm}

\smallskip

For $k\in{\mathcal C}_{\Delta}$, by (\ref{Delta(k)Blow})
 $$\langle{\rm Cm}(\Delta(k)),\,{\bf b}\rangle=
\langle{\rm Cm}(\Delta_0(\lambda,\tau)),\,{\bf
b}\rangle-\sum_{j=1}^nb_jk_j.$$
 If $(\Delta,{\bf b})$ is a mass linear
pair,
 by (\ref{Cmtruncated}) and Proposition \ref{PropLinear}, we have
$ \langle{\rm Cm}(\Delta_0(\lambda,\tau)),\,{\bf
b}\rangle=b_n\tau.$
 Thus,
$$\langle{\rm Cm}(\Delta(k)),\,{\bf
b}\rangle=\sum_{j=1}^{n+1}R_jk_j,$$
 where $R_j=b_n-b_j$, for $j=1,\dots, n$ and $R_{n+1}=b_n$. Hence, we
 have the following proposition:
  \begin{Prop}\label{Propsumrblow}
  Let $\Delta$ be the polytope obtained by truncating the standard $n$-simplex $S_n(\tau)$ by a horizontal hyperplane.
If $(\Delta,{\bf b})$ es a mass linear pair, then $\sum_jR_j=0$.
 \end{Prop}

\smallskip

{\it Remark.} When $n=3$ the   toric manifold $M_{\Delta}$ is
 $$M_{\Delta}=\big\{z\in{\Bbb
C}^5\,:\,|z_1|^2+|z_2|^2+|z_3|^2+|z_5|^2=\tau/\pi,\,\,|z_3|^2+|z_4|^2=\lambda/\pi
\big\}/{\Bbb T},$$
 where the action of ${\Bbb T}=(U(1))^2$ is
defined by
\begin{equation}\label{ActD}
 (a,b)(z_1,z_2,z_3,z_4,z_5)=(az_1,az_2,abz_3,bz_4,az_5),
 \end{equation}
  for
$a,b\in U(1)$.

We consider the following loops in the Hamiltonian group of
$(M_{\Delta},\,\omega_{\Delta})$
$$\psi_t[z]=[z_1e^{2\pi it},z_2,z_3,z_4,z_5],\;\;
\psi'_t[z]=[z_1,z_2e^{2\pi it},z_3,z_4,z_5],$$
$$\tilde\psi_t[z]=[z_1,z_2,z_3e^{2\pi it},z_4,z_5].$$
In \cite{V2} (Remark in Section 4), we gave formulas that relate
 the characteristic numbers
associated with these loops
$$I({\psi})=I({\psi'})=(-1/3)I({\tilde\psi}).$$
 So, for ${\bf b}=(b_1,b_2,b_3)\in{\mathbb Z}^3$,
 \begin{equation}\label{RemarkIpsi}
 I({\psi_{\bf b}})=(b_1+b_2-3b_3)I({\psi}).
\end{equation}
 By Proposition
\ref{PropLinear}, the vanishing of $I({\psi_{\bf b}})$ in
(\ref{RemarkIpsi}) is equivalent  to the fact that $(\Delta,\,{\bf
b})$ is a mass linear pair. This equivalence is a new checking of
  Theorem \ref{FinalThm}.

%%%%%%%%%%%%%%%%%%%%%%%%%%%%%%%%%%%%%%%%%%%%%%%%%%%%%%%%%%%%%%%%%%%%%%%%%%%%%%%%%%%%%%%%%%%%%%%%%%%%%%%%%%%%%%%%%%%%%%%%%%
%%%%%%%%%%%%%%%%%%%%%%%%%%%%%%%%%%%%%%%%%%%%%%%%%%%%%%%%%%%%%%%%%%%%%%%%%%%%%%%%%%%%%%%%%%%%%%%%%%%%%%%%%%%%%%%%%%%%%%%%%%

\end{document}